\def\bibart#1#2#3#4#5#6#7
{\bibitem{#1}
#2, #4, #5 #6 (#3), #7}

\def\bibcoll#1#2#3#4#5#6#7#8
{\bibitem{#1}
#2, #4, in: #5, #6, #7, #3, pp. #8}

\def\bibbook#1#2#3#4#5#6
{\bibitem{#1} #2, #4, #5, #6, #3}

\def\bibdiss#1#2#3#4#5#6
{\bibitem{#1}
#2 (#3), #4, #5, #6}

\def\qed{\hfill\ \rule{2mm}{2mm} }

\def\qex{\hfill\ \vbox{\hrule\hbox{\vrule\kern4pt\vbox{\kern4pt{}
\kern4pt}\kern4pt\vrule}\hrule}}

\documentclass[12pt]{article}

\usepackage{amssymb}

\newtheorem{guess}{Guess}[section]

\newtheorem{define}[guess]{Definition}
\newtheorem{prop}[guess]{Proposition}
\newtheorem{theorem}[guess]{Theorem}

\newtheorem{cor}[guess]{Corollary}

\newtheorem{lem}[guess]{Lemma}

\usepackage[pdftitle={FVP_FCP_for_Lexicographic_Sums}, pdfauthor={Bernd Schroeder}, baseurl={http://www.math.usm.edu/schroeder/}, colorlinks=true, pdfstartview=FitV, linkcolor=blue, citecolor=blue, urlcolor=blue, naturalnames=true, linktocpage=true, bookmarksnumbered=true, pdfpagelayout=SinglePage]{hyperref}

\begin{document}

\bibliographystyle{plain}

\title{Upper Covers of $2$-Chains
and of $2$-Antichains
in Sets of
Indecomposable Subsets}

\author{
\small Bernd S. W. Schr\"oder \\
\small School of Mathematics and Natural Sciences\\
\small The University of Southern Mississippi\\
\small
118 College Avenue, \#5043\\
\small Hattiesburg, MS 39406\\
}

\date{\small \today}

\maketitle

\begin{abstract}

We
prove that
there are arbitrarily large indecomposable ordered sets $T$
with a $2$-chain $C\subset T$ such that the smallest indecomposable
proper superset $U$ of $C$ in $T$ is $T$ itself.
Subsequently, we characterize all
such indecomposable ordered sets $T$ and
$2$-chains $C$.
We also prove the same type of result for $2$-antichains.

\end{abstract}


\section{Introduction}

In \cite{SchmTrot}, Schmerl and Trotter proved the following (in a more general context).

\begin{theorem}
\label{finitept3}

(See Theorem 2.2 in \cite{SchmTrot}.)
Let
$T$ be a finite indecomposable ordered set and let
$P\subset T$ be an indecomposable
ordered subset of $T$ with $4\leq |P|\leq |T|+2$.
Then there is an indecomposable ordered set
$U$ such that
$P\subset U\subseteq T$
and $|U|= |P|+2 $.
\qed

\end{theorem}

$2$-chains and $2$-antichains satisfy the definition of
indecomposability and, for them,
Theorem \ref{finitept3} fails.
Hence the requirement that $|P|\geq 4$ is needed.
This note characterizes all the ways in which
Theorem \ref{finitept3} fails for $2$-chains and for $2$-antichains.

\section{Upper Covers of Nonadjacent $2$-Chains}
\label{coverrelsec}

We first
consider certain
$2$-chains $C_2 =\{ a<b\} $ such that $a$ is not a lower cover of $b$.
Note that, for the family ${\cal X}$ below, the ordered sets
$P$ such that there are $a,b$ with $(P,a,b)\in {\cal X}$
include the ordered sets that are obtained as
finite convex indecomposable subsets
with at least $4$ elements of the infinite ordered sets in
\cite{brightlininf} as well as the ordered sets obtained from the
$3$-irreducible ordered sets $G_n $, $J_n $ and $H_n $
(see \cite{Trot}, p. 65) by deleting the elements
$c$ and $d$.

\begin{figure}

\centerline{
\unitlength 1mm 
\linethickness{0.4pt}
\ifx\plotpoint\undefined\newsavebox{\plotpoint}\fi 
\begin{picture}(90.5,36)(0,0)
\put(55,5){\circle*{1}}
\put(5,4.9){\circle*{1}}
\put(30,4.9){\circle*{1}}
\put(80,5){\circle*{1}}
\put(90,5){\circle*{1}}
\put(80,35){\circle*{1}}
\put(90,35){\circle*{1}}
\put(55,25){\circle*{1}}
\put(80,25){\circle*{1}}
\put(55,15){\circle*{1}}
\put(5,14.9){\circle*{1}}
\put(30,14.9){\circle*{1}}
\put(80,15){\circle*{1}}
\put(65,5){\circle*{1}}
\put(15,4.9){\circle*{1}}
\put(40,4.9){\circle*{1}}
\put(65,25){\circle*{1}}
\put(40,24.9){\circle*{1}}
\put(65,15){\circle*{1}}
\put(15,14.9){\circle*{1}}
\put(40,14.9){\circle*{1}}
\put(90,20){\circle*{1}}
\put(55,5){\line(0,1){10}}
\put(5,4.9){\line(0,1){10}}
\put(30,4.9){\line(0,1){10}}
\put(80,5){\line(0,1){10}}
\put(80,35){\line(0,-1){10}}
\put(55,15){\line(0,1){10}}
\put(80,15){\line(0,1){10}}
\put(55,15){\line(1,-1){10}}
\put(5,14.9){\line(1,-1){10}}
\put(30,14.9){\line(1,-1){10}}
\put(55,25){\line(1,-1){10}}
\put(65,5){\line(0,1){10}}
\put(15,4.9){\line(0,1){10}}
\put(40,4.9){\line(0,1){10}}
\put(65,15){\line(0,1){10}}
\put(40,14.9){\line(0,1){10}}
\put(65,4){\makebox(0,0)[t]{\footnotesize $a$}}
\put(15,3.9){\makebox(0,0)[t]{\footnotesize $a$}}
\put(40,3.9){\makebox(0,0)[t]{\footnotesize $a$}}
\put(55,26){\makebox(0,0)[cb]{\footnotesize $b$}}
\put(30,3.9){\makebox(0,0)[t]{\footnotesize $a'$}}
\put(80,4){\makebox(0,0)[ct]{\footnotesize $a$}}
\put(90,4){\makebox(0,0)[ct]{\footnotesize $\widetilde{a}$}}
\put(80,36){\makebox(0,0)[cb]{\footnotesize $b$}}
\put(90,36){\makebox(0,0)[cb]{\footnotesize $\widetilde{b} $}}
\put(40,25.9){\makebox(0,0)[cb]{\footnotesize $b$}}
\put(4,14.9){\makebox(0,0)[rc]{\footnotesize $b $}}
\put(29,14.9){\makebox(0,0)[rc]{\footnotesize $b' $}}
\put(41,14.9){\makebox(0,0)[lc]{\footnotesize $r $}}
\put(90,5){\line(-1,2){10}}
\put(90,35){\line(-1,-2){10}}
\put(90,35){\line(0,-1){15}}
\put(90,20){\line(0,-1){15}}
\put(80,35){\line(2,-3){10}}
\put(90,20){\line(-2,-3){10}}
\put(55,5){\line(1,2){10}}
\put(30,4.9){\line(1,2){10}}
\put(1,4.9){\makebox(0,0)[cc]{$N$}}
\put(51,4.9){\makebox(0,0)[cc]{$Y$}}
\put(26,4.9){\makebox(0,0)[cc]{$X $}}
\put(76,4.9){\makebox(0,0)[cc]{$Z$}}
\end{picture}
}

\caption{The first four elements $(P,a,b)$
in
the class ${\cal X}$ from Definition \ref{classX}.
We have $(N,a,b)\in {\cal X}$
by part \ref{classX1}
of Definition \ref{classX}
and, for $H\in \{ X,Y,Z\} $,
$(H,a,b)\in {\cal X}$ is obtained by applying
part \ref{classX2} or \ref{classX3}
of Definition \ref{classX}
to the set to the left of $H$.
The labeling of $X$ is used in the proof of Lemma \ref{2chnocovlem}.
}
\label{xtrafor2}

\end{figure}
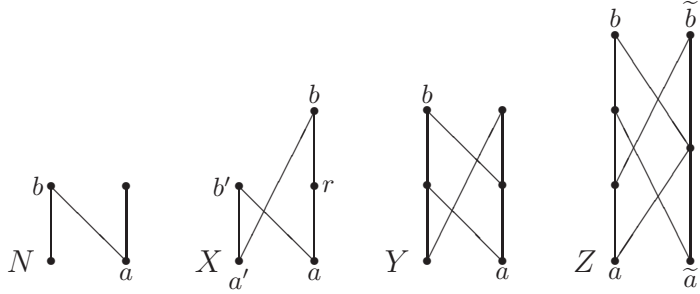

\begin{define}
\label{classX}

We define the family ${\cal X}$ of triples
$(P,a,b)$ of a
finite ordered set $P$, a minimal element $a\in P$
and a maximal element $b\in P$
by saying that
$(P,a,b)\in {\cal X}$
iff one of the following hold.
\begin{enumerate}
\item
\label{classX1}
$P$ is an $N$
and $a$ and $b$ are placed as in Figure \ref{xtrafor2}.
\item
\label{classX2}
There is a $(\widetilde{P} ,\widetilde{a} ,\widetilde{b} )\in {\cal X}$ such that
$P$ is obtained from $\widetilde{P}$ by attaching
$a$ as a new minimal element
below $\widetilde{P}\setminus \{ \widetilde{a}\} $
and
$b=\widetilde{b}$.
\item
\label{classX3}
There is a $(\widetilde{P},\widetilde{a},\widetilde{b})\in {\cal X}$ such that
$P$ is obtained from $\widetilde{P}$
by attaching $b$ as
a new maximal element
above $\widetilde{P}\setminus \{ \widetilde{b} \} $
and $a=\widetilde{a}$.

\end{enumerate}

\end{define}

We could immediately
show that, if $(P,a,b)\in {\cal X}$,
then the only indecomposable subset of
$P$ that contains $a$ and $b$ is $P$ itself.
However, a direct argument is
more technical than needed.
Hence we delay this discussion until after the proof of
Proposition \ref{2chfinal}.
We start by proving that certain $2$-chains
$\{ a,b\} $,
in which the elements are not covers of each other,
will be contained in
sets $P\subseteq T$ such that
$(P,a,b)\in {\cal X}$.
The hypothesis is a bit technical, but the overall
situation for $2$-chains
will be resolved in the proof of
Proposition \ref{2chfinal}.

\begin{lem}
\label{2chnocovlem}

Let
$T$ be a finite indecomposable ordered set with $|T|>2$ and let
$C_2 =\{ a<b\} \subset T$ be a chain with $2$ elements
such that $T\setminus \{ t\in T: a<t<b\} $ is series-decomposable.
Then
$C_2 $ is contained in a subset $H$ of $T$
such that
$(H,a,b)\in {\cal X}$ and $H$ is not isomorphic to $N$.

\end{lem}

{\bf Proof.}
Let
$R:=\{ t\in T: a<t<b\} $ and suppose, for a contradiction, that
$T$ is a finite indecomposable ordered set with $|T|>2$
such that the result does not hold and such that $R$ is as small as possible.

Consider the ordered set $S:=T\setminus R$.
By assumption, $S$ is series-decomposable.
Hence,
there are nonempty subsets $L, U\subset S$ such that
$S=L\oplus U$.
If $b\in L$, then $T=(L\cup R)\oplus U$, which is not possible.
If $a\in U$, then $T=L\oplus (R\cup U)$, which is not possible.
Thus $a\in L$ and $b\in U$.
Because $S$ contains no elements
strictly between $a$ and $b$, $a$ is maximal in $L$ and $b$ is minimal in $U$.
Moreover, because $T$ is indecomposable, we conclude that $R\not=\emptyset $.

Let $R_a :=\{ r\in R: r\not\geq L\} $
and let $R_b :=\{ r\in R: r\not\leq U\} $.
Note that both sets are nonempty, because otherwise
$T$ would be series-decomposable into $L\oplus (R\cup U)$
or into $(L\cup R)\oplus U$, respectively.
Pick $r_a \in R_a $
and $r_b \in R_b $
such that, if $R_a \cap R_b \not= \emptyset $, then
$r_a =r_b =:r$. Then there is an
$a''\in L$ such that $r_a \not\geq a'' $.
Let $a'\in L$ be a maximal element of $L$
such that $a'\geq a''$.
Then $a'$ is not comparable to $r_a $, because
$a'>r_a $ implies $a'>r_a>a$
(contradicting maximality of $a$ in $L$)
and
$a'<r_a $ implies $r_a >a'\geq a''$ (contradicting the choice of $a''$).
Similarly, there is a
$b'\in U$ that is not comparable to $r_b $ and minimal in $U$.

If $R_a \cap R_b \not= \emptyset $, then
$\{ a,a',r,b,b'\} $ is isomorphic to the ordered set
$X$ in Figure \ref{xtrafor2} and $(X,a,b)\in {\cal X}$, a contradiction
to the choice of $T$.
For the remainder, we can assume that
$R_a \cap R_b = \emptyset $.
That is, every element of $R$
is below $U$ or above $L$.
Therefore,
because $T=L\cup R\cup U$,
every element of $T$
is below $U$ or above $L$.

Now let
$\widetilde{a}\in L$ be maximal in $L$
and let $\widetilde{b}\in U$ be minimal in $U$, chosen so that
$\widetilde{R}:=\{ t\in T: \widetilde{a}<t<\widetilde{b}\} $
is as small as possible.
Note that $\widetilde{R}$ intersects neither $L$ nor $U$,
which means that
$\widetilde{R}\subseteq R$.
Because
$R':=\{ t\in T: a'<t<b'\} \subseteq R\setminus \{ r_a, r_b \} $,
we have that $|\widetilde{R}|\leq |R'|<|R|$, which means
$\widetilde{R}\subsetneq R$ and hence we have
$\widetilde{a}\not=a $ or $\widetilde{b}\not= b$.

Because every element of $T$
is below $U$ or above $L$,
every
element of $T$ is comparable to $\widetilde{a}$ or to $\widetilde{b} $.
Because $\widetilde{R}\subset R$, we have
$a<\widetilde{R}\setminus \{ \widetilde{a} \} $ and
$b>\widetilde{R}\setminus \{ \widetilde{b} \} $.
Let
$\widetilde{L}:=\{ t\in T\setminus \widetilde{R} :t<\widetilde{b} \} $
and let
$\widetilde{U}:=\{ t\in T\setminus \widetilde{R} :t>\widetilde{a} \} $.
Then
$\widetilde{a} $ is maximal in $\widetilde{L}$,
$\widetilde{b} $ is minimal in $\widetilde{U}$,
and
$\widetilde{L}\cap \widetilde{U}=\emptyset $.
Moreover, because every
element of $T$ is comparable to $\widetilde{a}$ or $\widetilde{b} $, we have
$\widetilde{L}\cup \widetilde{U}=T\setminus \widetilde{R} $.

Consider the case that
$\widetilde{a},\widetilde{b}$ can be chosen so that
$T\setminus \widetilde{R}$ is co-connected.
If
$\widetilde{L}<\widetilde{U}$, then we would have
$T\setminus \widetilde{R} = \widetilde{L} \oplus \widetilde{U}$, which cannot be.
Hence $\widetilde{L}\not<\widetilde{U}$, which means that there is
a maximal $\widetilde{\ell }\in \widetilde{L}$ that is
not comparable to a minimal $\widetilde{u}\in \widetilde{U} $.
If $\widetilde{\ell }
\in L$,
then
$\widetilde{u }
\not> L$, but
$\widetilde{u} \in \widetilde{U}$
also gives that $\widetilde{a}<\widetilde{u}$,
so that,
because $\widetilde{u}\not\in \widetilde{R} $,
$\widetilde{u}$ is not
smaller than
$\widetilde{b} $, implying that
$\widetilde{u}\not<U$, which contradicts the fact that $\widetilde{u} $
must be above $L$ or below $U$.
Similarly we exclude
$\widetilde{u }
\in U$.
Thus, $\widetilde{\ell} , \widetilde{u} \in R$, and,
in particular,
$\widetilde{\ell }
\not=
\widetilde{a}$ and
$\widetilde{u }
\not= \widetilde{b}$.
However, then
$\widetilde{N}:=\{ \widetilde{\ell }<\widetilde{b}>\widetilde{a}<\widetilde{u} \} $
is such that $(\widetilde{N}, \widetilde{a}, \widetilde{b}) \in {\cal X}$,
and subsequently
$(\widetilde{N} \cup \{ a,b\} , {a}, {b}) \in {\cal X}$,
contradicting the choice of $T$.

Therefore, $\widetilde{a},\widetilde{b}$
can only be chosen so that $T\setminus \widetilde{R}$ is series-decomposable.
Hence, by
choice of $T$ and because $|\widetilde{R}|<|R|$,
there is an ordered set $\widetilde{P}\subseteq T$ such that
$(\widetilde{P},\widetilde{a},\widetilde{b})
\in {\cal X}$.
With $P:=\widetilde{P}\cup \{ a,b\} $,
we have $(P,a,b)\in {\cal X}$,
independent of whether
$|P|=|\widetilde{P} |+1$
or
$|P|=|\widetilde{P} |+2$, a
final contradiction to the choice of $T$.
\qed

\begin{figure}

\centerline{
\unitlength 1mm 
\linethickness{0.4pt}
\ifx\plotpoint\undefined\newsavebox{\plotpoint}\fi 
\begin{picture}(146,66)(0,0)
\put(25,45){\circle*{1}}
\put(65,45){\circle*{1}}
\put(45,65){\circle*{1}}
\put(100,45){\circle*{1}}
\put(30,5){\circle*{1}}
\put(65,5){\circle*{1}}
\put(100,5){\circle*{1}}
\put(135,5){\circle*{1}}
\put(100,65){\circle*{1}}
\put(30,35){\circle*{1}}
\put(65,35){\circle*{1}}
\put(100,35){\circle*{1}}
\put(135,35){\circle*{1}}
\put(97.5,30){\circle*{1}}
\put(132.5,30){\circle*{1}}
\put(75,35){\circle*{1}}
\put(110,35){\circle*{1}}
\put(145,35){\circle*{1}}
\put(80,65){\circle*{1}}
\put(10,25){\circle*{1}}
\put(65,25){\circle*{1}}
\put(100,25){\circle*{1}}
\put(135,25){\circle*{1}}
\put(25,25){\circle*{1}}
\put(60,25){\circle*{1}}
\put(95,25){\circle*{1}}
\put(130,25){\circle*{1}}
\put(15,55){\circle*{1}}
\put(55,55){\circle*{1}}
\put(90,55){\circle*{1}}
\put(20,15){\circle*{1}}
\put(55,15){\circle*{1}}
\put(90,15){\circle*{1}}
\put(125,15){\circle*{1}}
\put(110,55){\circle*{1}}
\put(40,25){\circle*{1}}
\put(75,25){\circle*{1}}
\put(110,25){\circle*{1}}
\put(145,25){\circle*{1}}
\put(25,45){\line(-1,1){10}}
\put(65,45){\line(-1,1){10}}
\put(45,65){\line(1,-1){10}}
\put(100,45){\line(-1,1){10}}
\put(30,5){\line(-1,1){10}}
\put(65,5){\line(-1,1){10}}
\put(100,5){\line(-1,1){10}}
\put(135,5){\line(-1,1){10}}
\put(100,65){\line(-1,-1){10}}
\put(100,65){\line(1,-1){10}}
\put(30,35){\line(1,-1){10}}
\put(65,35){\line(1,-1){10}}
\put(100,35){\line(1,-1){10}}
\put(135,35){\line(1,-1){10}}
\put(80,65){\line(1,-1){10}}
\put(10,25){\line(1,-1){10}}
\put(65,25){\line(-1,-1){10}}
\put(100,25){\line(-1,-1){10}}
\put(135,25){\line(-1,-1){10}}
\put(25,45){\line(0,1){10}}
\put(15,55){\line(0,-1){10}}
\put(65,45){\line(0,1){10}}
\put(45,65){\line(0,-1){10}}
\put(100,65){\line(0,-1){10}}
\put(100,45){\line(0,1){10}}
\put(30,5){\line(0,1){10}}
\put(65,5){\line(0,1){10}}
\put(100,5){\line(0,1){10}}
\put(135,5){\line(0,1){10}}
\put(25,55){\circle*{1}}
\put(15,45){\circle*{1}}
\put(65,55){\circle*{1}}
\put(45,55){\circle*{1}}
\put(100,55){\circle*{1}}
\put(30,15){\circle*{1}}
\put(65,15){\circle*{1}}
\put(100,15){\circle*{1}}
\put(135,15){\circle*{1}}
\put(25,44){\makebox(0,0)[ct]{\footnotesize $a$}}
\put(65,44){\makebox(0,0)[ct]{\footnotesize $a$}}
\put(100,44){\makebox(0,0)[ct]{\footnotesize $a$}}
\put(30,4){\makebox(0,0)[ct]{\footnotesize $a$}}
\put(65,4){\makebox(0,0)[ct]{\footnotesize $a$}}
\put(100,4){\makebox(0,0)[ct]{\footnotesize $a$}}
\put(135,4){\makebox(0,0)[ct]{\footnotesize $a$}}
\put(14,55){\makebox(0,0)[rc]{\footnotesize $b$}}
\put(54,55){\makebox(0,0)[rc]{\footnotesize $b$}}
\put(89,55){\makebox(0,0)[rc]{\footnotesize $b$}}
\put(19,15){\makebox(0,0)[rc]{\footnotesize $b$}}
\put(54,15){\makebox(0,0)[rc]{\footnotesize $b$}}
\put(89,15){\makebox(0,0)[rc]{\footnotesize $b$}}
\put(124,15){\makebox(0,0)[rc]{\footnotesize $b$}}
\put(99,55){\makebox(0,0)[r]{\footnotesize $x$}}
\put(29,15){\makebox(0,0)[r]{\footnotesize $x$}}
\put(64,15){\makebox(0,0)[r]{\footnotesize $x$}}
\put(99,15){\makebox(0,0)[r]{\footnotesize $x$}}
\put(134,15){\makebox(0,0)[r]{\footnotesize $x$}}
\put(26,55){\makebox(0,0)[lc]{\footnotesize $x$}}
\put(66,55){\makebox(0,0)[lc]{\footnotesize $x$}}
\put(6,45){\makebox(0,0)[cc]{$N$}}
\put(40,45){\makebox(0,0)[cc]{$\widehat{N}$}}
\put(90,45){\makebox(0,0)[cc]{$B$}}
\put(20,5){\makebox(0,0)[cc]{$\widehat{B}$}}
\put(55,5){\makebox(0,0)[cc]{$\widetilde{B}$}}
\put(90,5){\makebox(0,0)[cc]{$B'$}}
\put(125,5){\makebox(0,0)[cc]{$B''$}}
\put(44,55){\makebox(0,0)[rc]{$\footnotesize \ell $}}
\put(14,45){\makebox(0,0)[rc]{$\footnotesize \ell $}}
\put(100,66){\makebox(0,0)[cb]{\footnotesize $w$}}
\put(30,36){\makebox(0,0)[cb]{\footnotesize $w$}}
\put(65,36){\makebox(0,0)[cb]{\footnotesize $w$}}
\put(100,36){\makebox(0,0)[cb]{\footnotesize $w$}}
\put(135,36){\makebox(0,0)[cb]{\footnotesize $w$}}
\put(96.5,30){\makebox(0,0)[rc]{\footnotesize $v_2 ^{<} $}}
\put(131.5,30){\makebox(0,0)[rc]{\footnotesize $v_2 ^{<} $}}
\put(45,66){\makebox(0,0)[cb]{\footnotesize $w$}}
\put(10,26){\makebox(0,0)[b]{\footnotesize $v_1$}}
\put(66,25){\makebox(0,0)[l]{\footnotesize $v_1$}}
\put(101,25){\makebox(0,0)[l]{\footnotesize $v_1$}}
\put(136,25){\makebox(0,0)[l]{\footnotesize $v_1$}}
\put(75,36){\makebox(0,0)[lb]{\footnotesize $v_2 ^{\not<} $}}
\put(110,36){\makebox(0,0)[lb]{\footnotesize $v_2 ^{\not<} $}}
\put(145,36){\makebox(0,0)[lb]{\footnotesize $v_2 ^{\not<} $}}
\put(111,55){\makebox(0,0)[l]{\footnotesize $\ell $}}
\put(41,25){\makebox(0,0)[l]{\footnotesize $\ell $}}
\put(76,25){\makebox(0,0)[l]{\footnotesize $\ell $}}
\put(111,25){\makebox(0,0)[l]{\footnotesize $\ell $}}
\put(146,25){\makebox(0,0)[l]{\footnotesize $\ell $}}
\put(80,66){\makebox(0,0)[b]{\footnotesize $u$}}
\put(24,25){\makebox(0,0)[rc]{\footnotesize $u$}}
\put(59,25){\makebox(0,0)[rc]{\footnotesize $u$}}
\put(94,25){\makebox(0,0)[rc]{\footnotesize $u$}}
\put(129,25){\makebox(0,0)[rc]{\footnotesize $u$}}
\put(10,25){\line(2,-1){20}}
\put(30,15){\line(0,1){10}}
\put(65,15){\line(0,1){10}}
\put(100,15){\line(0,1){10}}
\put(135,15){\line(0,1){10}}
\put(30,25){\line(0,1){10}}
\put(65,25){\line(0,1){10}}
\put(100,25){\line(0,1){10}}
\put(135,25){\line(0,1){10}}
\put(20,15){\line(1,2){5}}
\put(55,15){\line(1,2){5}}
\put(90,15){\line(1,2){5}}
\put(125,15){\line(1,2){5}}
\put(25,25){\line(1,2){5}}
\put(60,25){\line(1,2){5}}
\put(95,25){\line(1,2){5}}
\put(130,25){\line(1,2){5}}
\put(60,25){\line(3,2){15}}
\put(95,25){\line(3,2){15}}
\put(75,35){\line(-1,-2){10}}
\put(100,25){\line(1,1){10}}
\put(135,25){\line(1,1){10}}
\put(97.5,30){\line(1,-6){2.5}}
\put(132.5,30){\line(1,-6){2.5}}
\put(132.5,30){\line(5,2){12.5}}
\end{picture}
}

\caption{The forbidden sets in Lemma \ref{2chlem}, with one possible
placement of $b$ and $x$
indicated.
The only other possible placement of $b$ and $x$ is obtained by
$b$ taking the place of $x$
and vice versa and keeping all other points fixed.
The labelings in the sets
correspond to cases in the proof of Lemma \ref{2chlem}.
}
\label{2ch_forbid}

\end{figure}
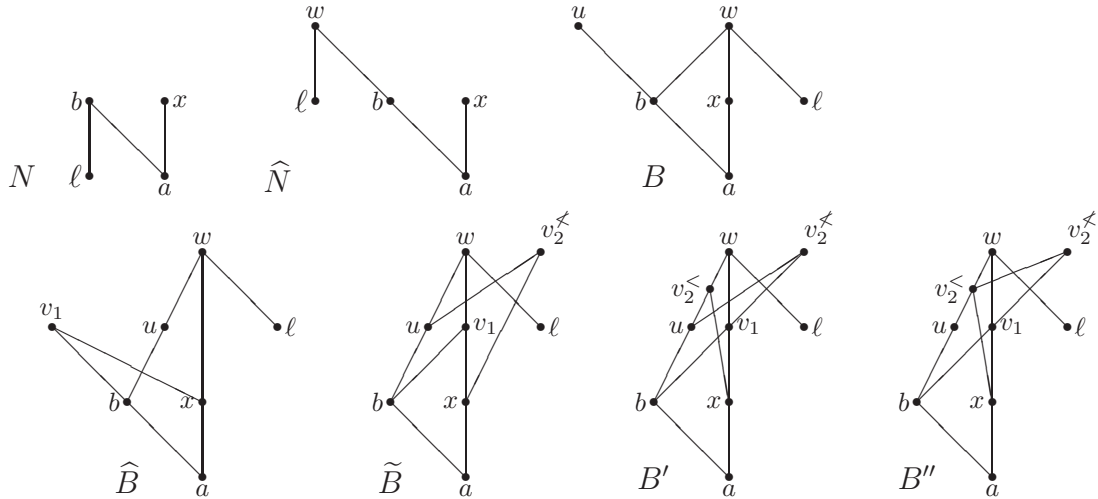

\section{Upper Covers of Adjacent $2$-Chains}

Now we prove that, if a chain in which both elements cover each other is contained in an
indecomposable ordered set, then the ordered set must contain one of the ordered sets in
Figure \ref{2ch_forbid} or its dual.

\begin{lem}
\label{2chlem}

Let
$T$ be a finite indecomposable ordered set with $|T|>2$ and let
$C_2 =\{ a<b\} \subset T$ be a chain with $2$ elements
such that $a$ is a lower cover of $b$.
Then
$C_2 $ is contained in a subset of $T$ that is isomorphic to one of the
ordered sets in Figure \ref{2ch_forbid} or to its dual.

\end{lem}

{\bf Proof.}
First note that,
if every strict upper bound of $a$ is an upper bound of $b$,
and
if every strict lower bound of $b$ is a lower bound of $a$,
then $C_2 $ is an order-autonomous subset of $T$, which is not possible.
Thus, without loss of generality, there is a
$p\in T$ such that $p>a$ and $p\not\geq b$.
(The dual case is not addressed in Figure \ref{2ch_forbid}, but, obviously,
runs along similar lines.)
Because $b$ is an upper cover of $a$, $b$ is not comparable to any element
$x$ with $a<x\leq p$.
In particular, $a$ has an upper cover $x\not= b$.
Let
$$C:=\{ x\in T: x {\rm \ is \ an \ upper \ cover \ of \ } a\} $$
and note that $|C|\geq 2$.

Consider the case that
there are
$x,y\in C$
that do not have the
same strict lower bounds.
Because we can, if needed, replace one of $x,y$ with $b$
and rename elements, we can assume that $y=b$.
If $b$ has a lower bound $\ell <b$ that is not a lower bound of
$x$, then
$\{ \ell <b>a<x\} $ is an ordered set $N$, see Figure \ref{2ch_forbid}.
If $x$ has a lower bound $\ell <x$ that is not a lower bound of
$b$, then
$\{ \ell <x>a<b\} $ is an ordered set
isomorphic to $N$, but, compared to
Figure \ref{2ch_forbid}, the positions of $b$ and $x$ are interchanged.
Thus, from here on, we can assume that
any two
$x,y\in C$
have the same strict lower bounds.

The argument above shows two characteristics of this proof.
First, we will continue to add hypotheses on the set $T$, typically indicated by
``from here on, we can assume."
These hypotheses typically are regarding the comparability of sets.
In Figure \ref{2ch_forbid}, elements of sets that will be introduced in the
future will be denoted with the corresponding lowercase letter and, possibly,
a superscript.
Second, although $b$ is an element of $C_2 $ and $x$ will be used to denote
a generic upper cover of $a$ that is not equal to $b$,
the roles of $x$ and $b$ in the sets depicted in
Figure \ref{2ch_forbid} will be interchangeable,
similar to the two versions of $N$ above.
This
interchangeability will be
indicated as needed.

Because any two
$x,y\in C$
have the same strict lower bounds
(and because any
element of $C$
is above $a$),
any $\ell \in T$ that is not comparable to either of $a$ or $b$
is not comparable to any
element of $C$.

Suppose, for a contradiction, that, for every
upper bound $w\in A:= \uparrow a\setminus \{ a\} $, every
lower bound $\ell <w$
is comparable to $a$ or $b$.
Let $p\in T\setminus A$ be such that $p$ is comparable to an element of $A$.
By definition of $A$, there is a $w\in A$ such that $p<w$.
By assumption,
$p$ is comparable to $a$ or $b$.
By definition of $A$, we have
$p\leq a$ or $p<b$.
In either case,
because any two elements of $C$ have the same strict lower bounds,
we infer that
$p<C$,
and hence
$p<A$.
We conclude that $A$ is nontrivially order-autonomous in $T$, a contradiction.
Thus there is an upper bound $w\in A= \uparrow a\setminus \{ a\} $ that has a
lower bound $\ell <w$ that is not comparable to either of $a$ or $b$.
Let
$$W:=\{ w\in T:w>a \wedge (\exists \ell <w) \ell \not\sim a\wedge \ell \not\sim b\}
\not=\emptyset .$$

Because any two
$x,y\in C$
have the same strict lower bounds,
we have that
$W\cap C=\emptyset $ and that any $\ell $ as in the
definition of $W$ is not comparable to
any element of $C$.
Also note that any element that is above an element of $W$
must be an element of $W$.

If there is
a $w\in W$
such that
$w\not>C$,
let $\ell <w$ be not comparable to either of $a$ and $b$
and choose $x\in C$ such that $w$ is comparable to one of
$b$ and $x$, but not the other.
Because
$\ell $ is not comparable to any
element of $C$,
we have that $\{ a,b\} $ is contained in an ordered set
isomorphic to $\widehat{N}$, see Figure \ref{2ch_forbid}
for the case $w>b$
(as indicated earlier, the case $w>x$ is obtained by switching $b$ and $x$).
Thus, from here on, we can assume that every
$w\in W$
is comparable to
all
elements of $C$, that is, $W>C$.

No two
$x,y\in C$
can have the same strict upper
bounds, because, otherwise,
$\{ x,y\} $ would be order-autonomous in $T$, which cannot be.
Hence, for any two
$x,y\in C$, there is a
$u\in T$ such that $u$ is a strict upper bound of one, but not the other.
Let
$$U:=\{ u\in T: (\exists x \in C)
[u>b \wedge u\not>x ]
\vee
[u>x \wedge u\not>b ]\} \not= \emptyset . $$

Note that, because $W>C$, we have that $U\cap W=\emptyset $.
If there is a $u\in U$ that is not
comparable to an element
$w\in W$,
let $x\in C$ be such that $u$ is above one of $b$ and $x$, but not the other.
Then
$\{ a,b\} $ is contained in an ordered set
isomorphic to $B$,
see Figure \ref{2ch_forbid} for the case $u>b$
(again, the case $u\not> b$ is obtained by switching $b$ and $x$).
Thus, from here on, we can assume that
$U<W$.
Next, we claim that
$$U=\{ u\in T: (\exists x,y \in C)
u>x \wedge u\not>y
\} . $$
The containment ``$\subseteq $" follows from the definition.
For the containment ``$\supseteq $,"
Let $t\in T$ and $x,y\in C$ be such that
$t>x$ and $t\not>y$.
In case $t>b$, we conclude that $t\in U$ because of the presence of $y$,
and,
in case $t\not>b$, we conclude that $t\in U$ because of the presence of $x$.
This proves the equality.

Let
$$V_1 :=\{ v_1 \in T\setminus W: v_1 >C, v_1 \not> U\} .$$
Note that no element of $U$ can be greater than or equal to any element of $V_1 $.

Suppose, for a contradiction, that $V_1 =\emptyset $.
Let $H:=C\cup U$ and let $p\in T\setminus H$ be comparable to an $h\in H$.
If $p<h$, then, by definition of $H$ and because $U\cap W=\emptyset $,
$c$ is below an upper cover
of $a$ and hence $p<C$ and then $p<H$.
If $p>h$, then $p$ is above an element of $C$.
Because $p\not\in U$, we obtain $p>C$,
and then,
because $V_1 =\emptyset $
and $W>U$,
we have $p>U$, which means that $p>H$.
We conclude that $H$ is nontrivially order-autonomous in $T$,
a contradiction.
Thus $V_1 \not=\emptyset $.
Note that, because $W>U$, no element of $V_1 $ is above any element of $W$.

If $V_1 \not<W $, then
there are a $v_1 \in V_1 $ and a $w\in W$ that are not comparable.
Because $v_1 \not>U $, there is a $u\in U$ that is not comparable to $v_1 $.
For this $u\in U$, there is an $x\in C$ such that $u$ is above one of $b$ and $x$,
but not the other.
This means that
$\{ a,b\} $ is contained in an ordered set
isomorphic to $\widehat{B}$, see Figure \ref{2ch_forbid} for the case
that $u>b$
(the case $u\not> b$ is obtained by switching $b$ and $x$).
Thus, from here on, we can assume that
$V_1 <W $.

Let
$$V_2 :=\{ v_2 \in T\setminus W: v_2 >U\} .$$
Note that no element of $V_2 $ is above any element of $W$ and
that no element of
$V_1 $ is above any element of $V_2 $.
Moreover, because $V_2 \cap U=\emptyset $, we have $V_2 >C$.

Suppose, for a contradiction, that $V_1 <V_2 $, and
consider the set $H:=C\cup U\cup V_1 $.
Let $p\in T\setminus H$.
If there is an $h\in H$ such that $p>h$, then
$p$ is above an element of $C$, so, because $p\not\in U$, we have
$p>C$, because
$p\not\in V_1 $ and $W>U$, we have $p>U$, and hence $p\in V_2\cup W>H$.
If there is an $h\in H$ such that $p<h$, then
(because $h\not\in W$)
$p$ is below an upper cover of $a$.
Hence $p<C$, which implies $p<H$ and $H$ is
nontrivially order-autonomous,
a contradiction.
Thus, $V_1 \not< V_2 $.

Let
\begin{eqnarray*}
V_2 ^{<}
& := &
\{ v_2 \in V_2 : v_2 <W\}
\\
V_2 ^{\not<}
& := &
V_2 \setminus V_2 ^{<}
\end{eqnarray*}

Note that no element of $V_2 ^{\not<} $ is below any element of $V_2 ^{<} $.

Now consider the case that
there is a $v_2 ^{\not<} \in V_2 ^{\not<} $ that is not an upper bound of
$V_1 $.
Then there is a $w\in W$ that is not above $v_2 ^{\not<} $,
there is a $v_1 \in V_1 $ that is not below $v_2 ^{\not<} $
(and hence not comparable to it)
and there is a
$u\in U$ that is not comparable to $v_1 $.
Consequently, using the comparabilities between the various sets that are already established,
we conclude that
$\{ a,b\} $ is contained in
an ordered set
isomorphic to
$\widetilde{B}$, see Figure \ref{2ch_forbid} for the case
that $u>b$
(the case $u\not> b$ is obtained by switching $b$ and $x$).
Thus, from here on, we can assume that
$V_1 <V_2 ^{\not<} $.
Because $V_1 \not<V_2 $, this means that
$V_1 \not< V_2 ^{<} $.

Suppose, for a contradiction, that
$V_2 ^{<} <V_2 ^{\not<} $, and
consider the set $H:=C\cup U\cup V_1 \cup V_2 ^{<} $.
Let $p\in T\setminus H$.
If there is an $h\in H$ such that $p>h$, then
$p$ is above an element of $C$, so,
because $p\not\in U$, we have $p>C$,
because $p\not\in V_1 $ and $W>U$, we have $p>U$,
so, because $p\not\in V_2 ^{<} $, we have
$p\in V_2 ^{\not<} \cup W>H$.
If there is an $h\in H$ such that $p<h$, then
(because $h\not\in W$)
$p$ is below an upper cover of $a$.
Hence $p<C$, which implies $p<H$ and $H$ is
nontrivially order-autonomous,
a contradiction.
Thus, $V_2 ^{<} \not< V_2 ^{\not<} $, which means that there are
$v_2 ^{<} \in V_2 ^{<} $ and $v_2 ^{\not<} \in V_2 ^{\not<} $ that are not comparable.
In particular, and this is all that will be used
in the following, this means that neither set is empty.

Let $v_2 ^{<} \in V_2 ^{<} $ be such that there is a $v_1 \in V_1 $
that is not below $v_2 ^{<} $
(and hence not comparable to it).
Then there is a
$u\in U$ that is not comparable to $v_1 $.

First, consider the case that
there is a $v_2 ^{\not<} \in V_2 ^{\not<} $
that is not comparable to $v_2 ^{<} $.
Then there is a $w\in W$ that is not above $v_2 ^{\not<} $.
Consequently, using the comparabilities between the various sets that are already established,
we conclude that
$\{ a,b\} $ is contained in an ordered set
$B'$, see Figure \ref{2ch_forbid} for the case
that $u>b$
(the case $u\not> b$ is obtained by switching $b$ and $x$).

Finally, if this is not the case, then
there is a $v_2 ^{\not<} \in V_2 ^{\not<} $ that is greater
than $v_2 ^{<}$.
Then there is a $w\in W$ that is not above $v_2 ^{\not<} $.
Consequently, using the comparabilities between the various sets that are already established,
we conclude that
$\{ a,b\} $ is contained in an ordered set
$B'$, see Figure \ref{2ch_forbid} for the case
that $u>b$
(the case $u\not> b$ is obtained by switching $b$ and $x$).
\qed

\vspace{.1in}

We can now summarize the situation for $2$-chains.

\begin{prop}
\label{2chfinal}

Let
$T$ be a finite indecomposable ordered set with $|T|>2$, let
$C_2 =\{ a<b\} \subset T$ be a chain with $2$ elements
and let $H$ be an indecomposable ordered subset of $T$ that
properly contains $C_2 $ such that there is no
indecomposable ordered subset
$U$ of $T$ with $C_2 \subsetneq U\subsetneq H$.
Then
$H$ is isomorphic to one of the
ordered sets in Figure \ref{2ch_forbid} or
$H$ is isomorphic to one of their duals,
or
$(H,a,b)\in {\cal X}$.

\end{prop}

{\bf Proof.}
Let
$S=H\setminus \{ t\in H: a<t<b\} $.
If $S$ is series-decomposable, the statement follows from
Lemma \ref{2chnocovlem}.
If $S$
is co-connected, then
$C_2 $ is not contained in a nontrivial order-autonomous subset of
$S$
(otherwise, $H$ would be decomposable).
We conclude that
$C_2 $ is contained in an indecomposable subset
$H'$ of $S$, and hence of $H$,
that is isomorphic to the index set of the canonical decomposition of $S$.
In particular, this means that $H'$ is not a $2$-chain
and, because $a$ is a lower cover of $b$ in $S$,
that $a$ is a lower cover of $b$ in $H'$.
Therefore, by
Lemma \ref{2chlem},
$C_2 $ is properly contained in an indecomposable subset of
$H'\subseteq H$ that is isomorphic to one of the
ordered sets in Figure \ref{2ch_forbid}
or to one of their duals.
The statement now follows from
the fact that this subset cannot be properly contained in $H$.
\qed

\vspace{.1in}

It can be checked that
none of the isomorphism types of the
ordered subsets
$H$ in Proposition \ref{2chfinal}
can be omitted:
One can prove that, for any two
chains $C=\{ a<b\} $ and $C'=\{ a'< b'\} $
and for ordered sets
$H$ (for $C$) and $H'$ (for $C'$)
as in Proposition \ref{2chfinal}, there is no
embedding of $H$ into $H'$ that maps $a$ to $a'$ and $b$ to $b'$.
Such an argument
is tedious, but can essentially be done ``by inspection."
Because of none of the isomorphism types of the
ordered subsets
$H$ in Proposition \ref{2chfinal}
can be omitted,
all bounds in
Corollary \ref{chcontainsummary} below are sharp.

\begin{cor}
\label{chcontainsummary}

Let $S$ be a finite ordered set
and consider the set ${\cal I}$ of all
indecomposable subsets of $S$, ordered by containment.
Let
$C\subset S$ ba a $2$-chain and
let $U$
be an upper cover of $C$ in ${\cal I}$.
If
$C=\{ a<b\} $ and $K$ is the longest chain from
$a$ to $b$, then $|U|\leq \max\{ 2|K| , 9\} $.

\end{cor}

{\bf Proof.}
The statement
follows from
Proposition \ref{2chfinal}.
\qed

\section{Upper Covers of $2$-Antichains in Sets of
Indecomposable Ordered Subsets}
\label{2accoverrelsec}

\begin{figure}

\centerline{
\unitlength 1.3mm 
\linethickness{0.4pt}
\ifx\plotpoint\undefined\newsavebox{\plotpoint}\fi 
\begin{picture}(136,25.5)(0,0)
\put(5,15){\circle*{1}}
\put(75,25){\circle*{1}}
\put(25,15){\circle*{1}}
\put(95,15){\circle*{1}}
\put(65,25){\circle*{1}}
\put(135,25){\circle*{1}}
\put(65,5){\circle*{1}}
\put(135,5){\circle*{1}}
\put(15,15){\circle*{1}}
\put(85,15){\circle*{1}}
\put(25,25){\circle*{1}}
\put(95,25){\circle*{1}}
\put(85,25){\circle*{1}}
\put(35,15){\circle*{1}}
\put(105,15){\circle*{1}}
\put(55,25){\circle*{1}}
\put(125,25){\circle*{1}}
\put(35,25){\circle*{1}}
\put(105,25){\circle*{1}}
\put(55,15){\circle*{1}}
\put(125,15){\circle*{1}}
\put(35,15){\line(0,1){10}}
\put(105,15){\line(0,1){10}}
\put(55,25){\line(0,-1){10}}
\put(125,25){\line(0,-1){10}}
\put(25,15){\line(0,1){10}}
\put(95,15){\line(0,1){10}}
\put(85,15){\line(0,1){10}}
\put(25,25){\line(-1,-1){10}}
\put(95,25){\line(-1,-1){10}}
\put(35,25){\line(-1,-1){10}}
\put(105,25){\line(-1,-1){10}}
\put(55,15){\line(1,1){10}}
\put(125,15){\line(1,1){10}}
\put(20,5){\circle*{1}}
\put(90,5){\circle*{1}}
\put(35,15){\line(-3,-2){15}}
\put(105,15){\line(-3,-2){15}}
\put(20,5){\line(1,2){5}}
\put(90,5){\line(1,2){5}}
\put(20,5){\line(-1,2){5}}
\put(90,5){\line(-1,2){5}}
\put(5,15){\line(3,-2){15}}
\put(20,3){\makebox(0,0)[ct]{\footnotesize $\ell $}}
\put(90,3){\makebox(0,0)[ct]{\footnotesize $\ell $}}
\put(4,15){\makebox(0,0)[rc]{\footnotesize $a $}}
\put(74,25){\makebox(0,0)[rc]{\footnotesize $a $}}
\put(14,15){\makebox(0,0)[rc]{\footnotesize $b $}}
\put(84,25){\makebox(0,0)[rc]{\footnotesize $b $}}
\put(35,15){\line(1,1){3.5}}
\put(105,15){\line(1,1){3.5}}
\put(55,25){\line(-1,-1){3.5}}
\put(125,25){\line(-1,-1){3.5}}
\put(65,25){\line(0,-1){20}}
\put(135,25){\line(0,-1){20}}
\put(45,20){\makebox(0,0)[cc]{$\cdots $}}
\put(115,20){\makebox(0,0)[cc]{$\cdots $}}
\multiput(20,5)(.0906735751,.0259067358){386}{\line(1,0){.0906735751}}
\multiput(90,5)(.0906735751,.0259067358){386}{\line(1,0){.0906735751}}
\put(66,25){\makebox(0,0)[lc]{\footnotesize $h$}}
\put(136,25){\makebox(0,0)[lc]{\footnotesize $h$}}
\put(66,5){\makebox(0,0)[lc]{\footnotesize $d$}}
\put(136,5){\makebox(0,0)[lc]{\footnotesize $d$}}
\put(90,5){\line(-3,4){15}}
\end{picture}
}

\caption{
Visualization of indecomposable V-covers, see Definition \ref{indecVcov}.
}
\label{2acd2}

\end{figure}
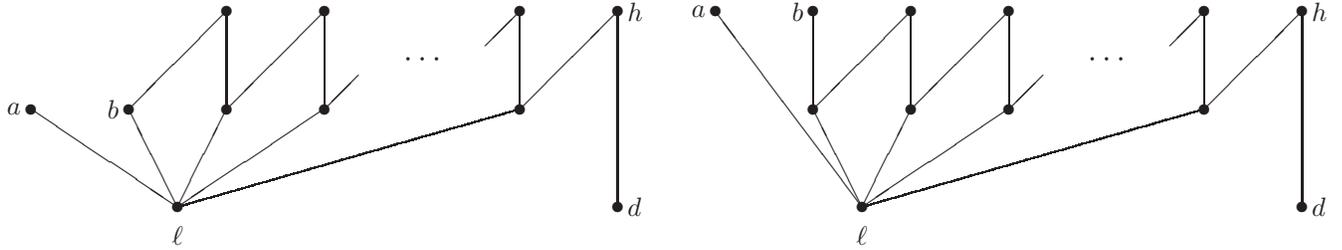

\begin{define}
\label{indecVcov}

An ordered set $C$ is called an {\bf indecomposable V-cover}
iff $C$ has exactly two minimal elements, $\ell $ and $d$,
the set $\uparrow \ell \setminus \{ \ell \} $
is composed of exactly two connected components, which are
a singleton $\{ a\} $ and a fence
$F$ from an element $b$ to an element $h$
such that $h$ is maximal in $F$
(with the equality $b=h$ being allowed);
and we have $\uparrow d\setminus \{ d\} =\{ h\} $.
The two types of indecomposable V-covers are depicted in
Figure \ref{2acd2}.

\end{define}

The name ``indecomposable V-cover"
comes from the following facts.
Clearly $a,b$ and $\ell $ form a fence that looks like a ``V."
Moreover, it is easy to check that
every indecomposable V-cover is indeed indecomposable.
Finally, we can easily see that, if $C$ is an indecomposable
V-cover, then the only indecomposable subset
$I\subseteq C$ that
properly contains $\{ a,b\} $ is $C$ itself:
Because such a subset $I$ must be connected, it must contain $\ell $.
Now, if any element of $F\setminus \{ b\} $ were not contained in $I$,
then, with
$B$ being the connected component of $F\cap I$ that contains $b$, we
would have that
$\{ a\} \cup B$ is nontrivially order-autonomous in $I$, which cannot be.
Hence $\{ a,\ell \} \cup F\subseteq I$, and, because
$\{ a,\ell \} \cup F$ is series-decomposable, we must have
$d\in I$ and hence $C=I$.

\begin{prop}
\label{2aclem}

Let
$T$ be a finite indecomposable ordered set with $|T|>2$ and let
$A_2 =\{
a,b
\} \subset T$ be an antichain with $2$ elements.
If $d(
a,b
)>2$, then any
smallest indecomposable ordered subset
$I$ of $T$ that contains $A_2 $
is a
fence from
$a$ to
$b$
with $d(
a,b
)+1$ elements.
If $d(
a,b
)=2$, then any
smallest indecomposable ordered subset
$I$ of $T$ that contains $A_2 $
is
either a fence with at least $4$ elements, or it is
isomorphic or dually isomorphic
to an indecomposable V-cover in which $a$ and $b$ are as in
Definition \ref{indecVcov}
(see Figure \ref{2acd2}) or in which the roles of $a$ and $b$ are
interchanged.

\end{prop}

{\bf Proof.}
The case $d(a,b)>2$ is already discussed at the start of this section.
We are left to
consider the case $d(a,b)=2$.

Let
$I$
be a smallest indecomposable ordered subset
of $T$ that contains $A_2 $.
If $I$ does not contain a common upper or lower bound of $a$ and $b$, then
an argument similar to the argument at the start of this section
shows that $I$ is a fence.
(Surprisingly, this case can occur: Consider $N=\{ a<f_2 >f_3 <b\} $ with an additional
element $\ell <a,b$ added.)
This leaves the case that $I$ contains a common upper or lower bound of $a$ and $b$.
Without loss of generality, assume that
$I$ contains a common lower bound of $a$ and $b$.
(The other case is handled with the dual argument.)

%
%
%

Let $L:=\{ x\in I: x <a,b\} \not= \emptyset $,
let $U:=\{ x\in I: x>a,b\} $
and let $H:=\{ x\in I: L<x<U\} \supseteq \{ a,b\} $ be the set of all
elements between $L$ and $U$ in $I$.
Let $A$ be the connected component of $H$ that contains $a$ and
let $B$ be the connected component of $H$ that contains $b$.
Note that $A$ could be equal to $B$.
If $|B|>1$, then, because $B$ cannot be order-autonomous in $I$, there
must be an element in $B$ that has a strict upper or lower bound
in $I$ that is not in
$L\cup B\cup U$.
Similarly, if $|A|>1$,
there is an element of $A$ that has a strict
upper or lower
bound in $I$ that is not in
$L\cup A\cup U$.
Finally, in case $|A|=|B|=1$, because $\{ a,b\} $ cannot be
order-autonomous in $I$,
$a$ must have
a strict upper or lower bound in $I$ that is not in $L\cup A\cup U$,
or $b$ must have
a strict upper or lower bound in $I$ that is not in $L\cup B\cup U$.
Note that all these upper and lower bounds are not in $L\cup H\cup U$.
Moreover, note that
the existence of a strict upper bound $s$ for an element of $A$ (or $B$)
such that
$s\not\in L\cup H\cup U$ implies $U\not= \emptyset $,
because otherwise $s$ would be in $A$ (or $B$).
Because
we can switch the roles of $a$ and $b$,
and because we can work with the dual ordered set if needed,
without loss of generality, we can assume that
$B$
contains an
element that has a strict lower bound in $I$
that is not in $L\cup B\cup U$.

Let $s$ be the shortest distance in $B$ from the element
$b$ to
an
element of $B$ that has a strict lower bound in $I$
that is not in $L\cup B\cup U$.
In case $a\in B$, because the roles of $a$ and $b$ can be switched,
we can assume that the distance in $B$ from the element $a$ to
any
element of $B$ that has a strict lower bound in $I$
that is not in $L\cup B\cup U$ is at least $s$.
Let $F\subseteq B$ be a
fence of length $s$ in $B$ that goes from
$b$ to an element $h\in B$ that has a lower bound
$d\in I\setminus (L\cup B\cup U)$.
Because
$a$ and $b$ do not have common upper or lower bounds in $B$,
$F$ has length $s$ and $F$ is contained in $B$,
we conclude that
$a$ is not comparable to any element of $F\cup \{ d\} $
and that
$d$ is not comparable to any element of $F\setminus \{ h\} $.
Because $d\not\in L$, no element of $L$ is above $d$.
Because $d<h$, we have $d<U$.
Because $d<h\in B$, $d\not\in B\subseteq H$ and $B$ is a connected
component of $H$, the point $d$ cannot be in $H$.
Because $d<U$, and $d\not\in H$, the point
$d$ cannot be above all elements of $L$.
Therefore, there is an
$\ell \in L$ that is incomparable to $d$.
This means that $C:=\{ a,b,d,\ell \} \cup F\subseteq I$
is an indecomposable V-cover that contains $a$ and $b$.
Because $I$ is smallest possible, we obtain that $I=C$.
\qed

\end{document}